\documentclass[11pt]{article}

\newcommand{\mP}{\mathbb{P}}
\newcommand{\mE}{\mathbb{E}}
\newcommand{\mF}{\mathcal{F}}
\newcommand{\wt}{\widehat}
\newcommand{\md}{\mathcal{P}}
\newcommand{\sd}{\md = \left \{ F : \int_{0}^{\infty} x \,\, dF(x) \leq 1 \right\}}
\newcommand{\sdeps}{\text{closure} \left[ \left \{ F : \int_{0}^{\infty} x \,\, dF(x) \leq 1 + \frac{1}{n} \right \} \right] }
\newcommand{\af}{\Phi}
\newcommand{\dm}{\delta_{\alpha} p + \delta_{\beta} (1 - p) = F}
\newcommand{\dmE}{\alpha p + \beta (1 - p) \leq 1}
\newcommand{\closedsetcap}{\bigcap_{n=1}^{\infty} \sdeps }
\newcommand{\compK}{K_{\varepsilon}}

\usepackage[utf8]{inputenc}
\usepackage{csquotes}
\usepackage{amsmath, amsthm, amssymb, amsfonts}
\usepackage{geometry}
\usepackage{leftindex}
\usepackage{hyperref}
\usepackage{graphicx}
\usepackage{cite}
\usepackage{mathrsfs}
\usepackage{tikz}
\usepackage{enumitem}

\newtheorem{theorem}{Theorem}[section]
\newtheorem{lemma}[theorem]{Lemma}

\theoremstyle{definition}

\title{An Equivalent Conjecture To Feige's Conjecture}
\author{Metin Dürr \footnote{My legal name is Bozkurt, however, due to the political associations with this name I prefer to be referenced by the surname of my future wife, i.e., Dürr.}\, \footnote{metin.duerr.research@proton.me}}
\date{\today}

\begin{document}

\maketitle
\begin{abstract}
\noindent
Let $X_{1}, ..., X_{n}$ be arbitrary non-negative independent random variables with respective expected values $\mu_{i}$ at most one. We sketch but do not prove an equivalent conjecture to Feige's Conjecture
\[
\mP \left( \sum_{i=1}^{n} X_{i} < \mu + 1 \right) \geq \exp (-1),
\]
where $\mu$ is the expected value of the sum of the random variables. We show by a simple example how this inequality finds use in mathematical finance.
\end{abstract}

\section{Introduction}
Many probabilistic inequalities, e.g., Markov's inequality, bound the probability mass above a given constant of a given random variable $X$. Moreover, additional conditions might be imposed, e.g., the existence of the variance of $X$. In contrast, Feige's Conjecture gives a lower bound on the probability mass of $X$ below a certain threshold and does not assume that the higher moments of $X$ exist. \newline
This research was motivated by a conjecture of Feige \cite{Feige}. For arbitrary non-negative independent random variables $X_{1}, ..., X_{n}$ with the respective expected values $\mu_{i}$ at most one the conjecture states
\begin{equation}
\mP \left( \sum_{i=1}^{n} X_{i} < \mu + 1 \right) \geq \lambda,
\end{equation}
where $\lambda = \exp(-1)$. In \cite{Feige} Feige proved a bound of approximately $\lambda = 0.0769$ by applying an algorithm, which reduces the support of a random variable to a discrete set with two elements, thereby enabling a case analysis. This bound was further improved by Garnett \cite{Garnett} to $\lambda = 0.14$ by considering the first four moments. The best known bound the author was able to find is $\lambda = 0.1798$ due to Guo et al. \cite{Guo}, who use an optimization approach paired with the Berry-Essen Theorem. \newline
Moreover, we want to note that Feige's Conjecture was proven true for discrete log-concave distributions by Alqasem et al. \cite{Al} and for identically distributed random variables by Egozcue et al. \cite{Ego}.
In this paper we give a different argument than the one found in~\cite{Feige} why it suffices to study discrete two-valued random variables, namely we utilize Prokhorov's Theorem. Moreover, we sketch an equivalent conjecture to Feige's Conjecture.

\section{Notation}
Throughout this paper we denote by large latin letters like $X$ arbitrary non-negative independent random variables and write r.v. as a shorthand. Moreover, we denote by $\mu$ expected values. We always assume that the random variables are defined on a common probability space $(\Omega, \mF, \mP)$.

\section{An Equivalence To Feige's Conjecture} \label{sb}
We first show that studying discrete two valued random variables suffices. We go through our argument in dimension one to avoid cumbersome notation. The ideas generalise to any dimension. Consider the set
\begin{equation}
\sd,
\end{equation}
where $F$ is the cumulative distribution function of a r.v. $X$. We equip $\md$ with the weak* topology. This set is convex and the extreme points are Dirac measures, such that
\begin{equation}
\dm,
\end{equation}
and
\begin{equation}
\dmE,
\end{equation}
which includes singletons. We refer to \cite{Winkler} for a proof of this. We will show that this set is also compact. First, note that $\md$ is closed since
\begin{equation}
\md = \closedsetcap.
\end{equation}
Now let $\varepsilon > 0$ and $\compK$ be a compact subset of $[0, \infty)$, such that for $a > \frac{\mu}{\varepsilon}$ we have \newline $A := [0, a) \subset \compK$. By Markov's Inequality we see
\begin{equation}
F(\compK) \geq \mP (X < a) \geq 1 - \frac{\mu}{a} > 1 - \varepsilon.
\end{equation}
Thus, $\md$ is tight. Moreover, since $[0, \infty)$ is a Polish space so is $\md$ and we see by Prokhorov's Theorem that $\md$ is compact. \newline
Consider now the function
\begin{equation}
\begin{split}
\af : \,\, &\md \rightarrow [0, 1], \\
&F_{X} \mapsto H_{X}(1 + \delta) := \mP (X \geq 1 + \delta).
\end{split}
\end{equation}
This function is affine since for $\lambda \in [0, 1]$ and $F_{X}, F_{Y} \in \md$ we have
\begin{equation}
\af (\lambda F_{X} + (1 - \lambda) F_{Y}) = (\lambda H_{X} + (1 - \lambda) H_{Y})(1 + \delta) = \lambda H_{X}(1 + \delta) + (1 - \lambda) H_{Y}(1 + \delta).
\end{equation}
Since for any $r \in \mathbb{R}$ the set
\begin{equation}
\left \{ t \in \mathbb{R} : \mP (X \geq t) < r \right \}
\end{equation}
is open in $\mathbb{R}$ we see that $H_{X}(1 + \delta)$ is upper semi-continuous and so is $\af$. By Bauer's Maximum Principle we know that $\af$ attains its maximum at the extreme points of $\md$ and thus $\mP (X < 1 + \delta)$ attains its minimum at the extreme points of $\md$. The minimum candidates for any dimension $n$ are
\begin{equation}
X_{i} = \begin{cases}
0, &p_{i}, \\
\frac{1}{1 - p_{i}}, &1 - p_{i}.
\end{cases}
\end{equation}
Note that we may restrict to the case $\mE (X_{i}) = 1$, since we can always define variables $Z_{i} := X_{i} + \wt{\mu_{i}}$ such that $\mu_{i} + \wt{\mu_{i}} = 1$ and we have
\begin{equation} \label{probeq}
\mP \left( \sum_{i=1}^{n} X_{i} < \mu_{1} + ... + \mu_{n} + \delta \right) = \mP \left( \sum_{i=1}^{n} Z_{i} < n + \delta \right).
\end{equation}
Intuitively, we assume that to minimize the probability~\ref{probeq} all the mass should be above $n + \delta$. This, however, is not the case for $\delta$ small enough, e.g., $\delta = 0.1$. The next Lemma gives us a necessary condition for our intuition to hold true, namely $\delta \geq 1$. We omit the proof, since it is a straightforward application of the AM-GM inequality.
\begin{lemma} \label{prob}
Let $0 < p_{1}, ..., p_{n} < 1$ and $n \in \mathbb{N}$. There exists no solution to the problem
\begin{enumerate}
\item $1 - \prod_{i=1}^{n} (1 - p_{i}) < \left( \frac{n}{n + 1} \right)^{n}$
\item $\sum_{i=1}^{n} \frac{1}{1 - p_{i}} \geq n + 1$
\end{enumerate}
\end{lemma}
There might, however, exist a minimizing mass distribution, where not all the mass is to the right of $n + 1$. To show this, one has to prove that all possible optimization problems for $n \in \mathbb{N}$ are bounded by $\exp(-1)$. We show what we mean for the case $n = 2$. We have the following two problems to consider.
\begin{enumerate}
\item $\inf p_{1}p_{2} + p_{1}(1-p_{2}) = p_{1}$
\item $p_{1} \geq \frac{2}{3}$
\end{enumerate}
and
\begin{enumerate}
\item $\inf p_{1}p_{2} + p_{1}(1-p_{2}) + (1 - p_{1})p_{2} = 1 - (1 - p_{1})(1 - p_{2})$
\item $\frac{1}{1 - p_{1}} + \frac{1}{1 - p_{2}} \geq 3$
\item $0 < p_{1}, p_{2} < \frac{2}{3}.$
\end{enumerate}
The solution to both of these problems is greater than $\exp(-1)$, where the second problem is already the statement of Lemma~\ref{prob}. If one succeeds in proving this for general $n$, one has succeded in proving Feige's Conjecture.

\section{An Example In Financial Mathematics}
Suppose we have invested into $n$ stocks of type $A$ at time $t_{0} = 0$. At some future time $t_{1} > t_{0}$ the value of each individual stock is higher than at time $t_{0}$. From the time $t_{1}$ onward let $X_{i}(t)$ describe the profit we would gain if we sold the stock $i$ at time $t \geq t_{1}$, where we immediately sell the stock once the profit drops to zero, i.e., the random variables $X_{i}(t)$ will never be negative, excluding transaction costs. We know from historical data, that a portfolio of type $A$ stocks has an expected value of $\mu$ units of currency per invested currency for the time $T > t_{1}$. We may ask ourselves, what the probability is that our portfolio will outperform this expected value $\mu$ by an amount of $\delta \geq 1$ at time $T$. This probability is given by
\begin{equation} \label{finineq}
1 - \exp ( - 1 ) \geq \mP \left( \sum_{i=1}^{n} X_{i}(T) \geq \mu + \delta \right).
\end{equation}
Of course, the historical data takes into account that the portfolio of type $A$ stocks can be negative while our strategy does not. Nonetheless, we may choose a suitable $\delta$ to counterbalance this and interpret the inequality~\ref{finineq} as a bound on the probability of outperforming the expected value $\tilde{\mu} = \mu + \delta$ of our strategy.

\section{Conclusion}
We have given a topological argument why considering discrete two-valued random variables is enough to prove Feige's Conjecture. Moreover, we sketched an equivalent Conjecture to Feige's Conjecture.

\bibliographystyle{plain}
\bibliography{references}

\begin{thebibliography}{1}

\bibitem{Al}
A.~Alqasem, H.~Aravinda, A.~Marsiglietti, and J.~Melbourne.
\newblock On a \uppercase{C}onjecture of \uppercase{F}eige for
  \uppercase{D}iscrete \uppercase{L}og-\uppercase{C}oncave
  \uppercase{D}istributions.
\newblock {\em SIAM Journal on Discrete Mathematics}, 38(1):93--102, 2024.

\bibitem{Ego}
M.~Egozcue and L.~F. García.
\newblock A simple proof of \uppercase{F}eige's conjecture for identically
  distributed random variables.
\newblock {\em https://ssrn.com/abstract=5188125}, 2025.

\bibitem{Feige}
U.~Feige.
\newblock On \uppercase{S}ums of \uppercase{I}ndependent \uppercase{R}andom
  \uppercase{V}ariables with \uppercase{U}nbounded \uppercase{V}ariance, and
  \uppercase{E}stimating the \uppercase{A}verage \uppercase{D}egree in a
  \uppercase{G}raph.
\newblock {\em SIAM Journal on Computing}, 35(4):964--984, 2006.

\bibitem{Garnett}
B.~Garnett.
\newblock Small \uppercase{D}eviations of \uppercase{S}ums of
  \uppercase{I}ndependent \uppercase{R}andom \uppercase{V}ariables.
\newblock {\em Journal of Combinatorial Theory, Series A}, 169(3), 2018.

\bibitem{Guo}
J.~Guo, S.~He, Z.~Ling, and Y.~Liu.
\newblock Bounding probability of small deviation on sum of independent random
  variables: \uppercase{C}ombination of moment approach and
  \uppercase{B}erry-\uppercase{E}ssen theorem, 2020.

\bibitem{Winkler}
G.~Winkler.
\newblock \uppercase{Extreme Points Of Moment Sets}.
\newblock {\em Mathematics of Operations Research}, 13(4):581--587, 1988.

\end{thebibliography}
\end{document}